\documentclass[12pt]{article}

\title{Stochastic differential equations
of second order
with a small parameter.
\thanks{
The research is funded by the grant of the Government of Russian Federation $n^\circ$14.A12.31.0007 
and  by The National Research  Tomsk State University
Academic D.I. Mendeleev Fund Program (NU 8.1.55.2015 L) in 2014- 2015.
}
}

\author{Kamenskii, M.,
\thanks{
 Voronezh State University,
Universitetskay pl. 1,
 394063 Voronezh, Russia,
 e-mail: Mikhailkamenski@mail.ru}
Quincampoix, M.
\thanks{
Laboratoire de Math\'ematiques de Bretagne Atlantique,
 CNRS UMR 6205,
Universit\'e de Bretagne Occidentale,
6 avenu Victor Le Gorgeu - CS 93837
29238 BREST Cedex 3, FRANCE,
e-mail : Marc.Quincampoix@univ-brest.fr
}
\and
S. Pergamenshchikov\thanks{
 Laboratoire de Math\'ematiques Raphael Salem,
CNRS UMR 6085,
 Avenue de l'Universit\'e, BP. 12,
  Universit\'e de Rouen,
   F76801, Saint Etienne du Rouvray, Cedex France,
and 
National Research University - Higher School of Economics, Laboratory of Quantitative Finance,
Moscow, Russia,
 e-mail: Serge.Pergamenchtchikov@univ-rouen.fr}
}

\usepackage{srcltx}
\usepackage{amssymb}
\usepackage{amsfonts}
\usepackage{amsmath}
\usepackage{amsthm}
\usepackage{enumerate}
\usepackage{multicol}

\newtheorem{theorem}{Theorem}[section]
\newtheorem{proposition}[theorem]{Proposition}
\newtheorem{lemma}[theorem]{Lemma}

\newtheorem{remark}{Remark}[section]

\newcommand\cC{{\cal C}}

\newcommand\ve{\varepsilon}

\def\bbr{{\mathbb R}}

\def\text#1{\hbox{#1}}
\def\proof{{\noindent \bf Proof. }}
\def\endproof{\mbox{\ $\qed$}}

\def\E{{\bf E}}
\def\P{{\bf P}}
\def\C{{\bf C}}

\def\H{{\bf H}}

\def\a{{\bf a}}
\def\b{{\bf b}}
\def\B{{\bf B}}

\def\AP{{\bf AP}}

\newcommand{\wt}{\widetilde}

\def\Chi{{\bf 1}}
\def\d{\mathrm{d}}
\def\build #1_#2{\mathrel{\mathop{\kern 0pt #1}\limits_{#2}}}
\newcommand{\zs}[1]{{\mathchoice{#1}{#1}{\lower.25ex\hbox{$\scriptstyle#1$}}
{\lower0.25ex\hbox{$\scriptscriptstyle#1$}}}}
\numberwithin{equation}{section}

\begin{document}

\maketitle

\begin{abstract}
We consider  boundary value problems
for stochastic differential equations of second order with a small parameter. 
For this case we prove a special existence and unicity theorem  for strong
solutions. The asymptotic behavior of these solutions as small parameter goes to zero
is studied. The stochastic averaging theorem for such equations is shown.
The limits in the explicit form for the solutions as a small parameter goes to zero
are found.
\end{abstract}
\vspace*{5mm}

\vspace*{5mm}

\bibliographystyle{plain}
\renewcommand{\columnseprule}{.1pt}

\newpage

\section{Introduction}\label{sec:1}

This paper is devoted to a limiting problem for the stochastic differential equations
with the small parameters. Such problems remount to the papers on the asymptotic theory
for the ordinary differential equations. There exist two general methods for study
these problems: the Bogolubov averaging principle \cite{BoMi}
 and the singularly perturbation Tikhonov method \cite{Ti}-\cite{Ti1}.
In the first case the quick variables have not limit, but in the second case the quick variables
go to some limit. Usually, to find this limit one has to replace in the equations the small parameter
by zero, and the solution of the limiting problem gives the limit for the quick variables.
 The situation becomes essentially difficult if we consider the boundary problems. 
So, let us consider, for example, the ordinary second order differential equation 
with the small parameter on the interval $[0,1]$, i.e.
\begin{equation}\label{sec:1.0} 
\ddot{x}^{\ve}(t)=\varepsilon\left( px^{\varepsilon}(t)+f(t)\right)\,,
\end{equation}
where $p>0$.
In this paper we consider the following boundary conditions:

\noindent {\sl First Neumann boundary conditions}

\begin{equation}\label{sec:2.2}
 x(0)=x(1)=0\,.
\end{equation}

\noindent {\sl Second  Neumann boundary conditions}

\begin{equation}\label{sec:2.3}
 \dot{x}(0)=\dot{x}(1)=0\,.
\end{equation}

\noindent {\sl Periodic boundary conditions}

\begin{equation}\label{sec:2.4}
x(0)=x(1)
 \quad\mbox{and}
\quad
\dot{x}(0)=
\dot{x}(1)\,.
\end{equation}
Note that if we replace in the equation \eqref{sec:1.0} the parameter
$\varepsilon$ by zero, we obtain the following limit problem 
$$
\ddot{x}=0\,.
$$
It should be noted that this equation with the first boundary condition has
unique solution $x=0$. 	However, for the boundary conditions 
\eqref{sec:2.3}
or \eqref{sec:2.4} this  limiting problem 
has not unique solution, since any constant satisfies this equation,
 i.e. in this case the limiting problem 
does not give any information about the asymptotic (as $\ve\to 0$)
 behavior of the solution of the equation \eqref{sec:1.0}.
 Therefore, a natural question arises, haw to calculate
this constant. For the boundary conditions \eqref{sec:2.3}
or \eqref{sec:2.4}  
through the Green functions method (see Remark~\ref{Re.sec:Av.0} in 
Section~\ref{sec:Av})  one can show, 
that for any $0< t<1$
\begin{equation}\label{sec:2.4-1}
\lim_\zs{\varepsilon\to 0}\,x^{\ve}(t)=
-\frac{1}{p}\int^{1}_\zs{0} f(s)\d s \,.
\end{equation}

\noindent There exists a vast literature on the stochastic differential equations with small
parameters and its applications (see, for example, \cite{FrWe}, \cite{KaPe}
and \cite{Ku} for details and references therein).
It turns out that the small parameter methods for stochastic differential equations 
are useful
for the optimal stochastic control problems (\cite{AlBa}, \cite{ArGa}, \cite{Be}),
for stochastic volatility financial markets 
(\cite{AiMyZh}, \cite{AiMyZh1}, \cite{FoPaSi}), 
for the statistical estimation in  multi-scale stochastic systems
(\cite{PaPaSt}, \cite{PaSt} \cite{PaSt1}) and 
for many other applied problems.
Usually, one considers initial or terminal conditions  for the stochastic differential
 equations.
In this paper we propose to extent both these problems, i.e. we introduce
the stochastic differential equation of the second order with the boundary conditions
of the forms \eqref{sec:2.2} -- \eqref{sec:2.4}. For such models we study 
the stochastic averaging principle. This problem is well known for
the diffusion processes (for example, \cite{FrWe} or \cite{KhaKr}).
To understand the averaging principle problem we start with the following example.
We consider 
 on the interval $[0,\ve]$ (for some small parameter $0<\ve<1$) the linear
stochastic differential
equation of second order with
 quickly oscillating
and ``small'' random force (``white noise''), i.e.
$$
\ddot{y}^{\ve}(u)
=
p y^{\ve}(u)
+
f(u/\ve)
+
\delta(u/\ve) \sqrt{\ve} \dot{w}_\zs{u}\,,
$$
where $p>0$ is some fixed parameter 
 and $\delta$ is a $\bbr\to\bbr$ square integrated function, i.e.
\begin{equation}\label{sec:2.4-1-0}
\|\delta\|^{2}=
\int^{1}_\zs{0}\,\delta^{2}(t)\d t\,<\,\infty\,.
\end{equation}
Moreover,  
 $(w_\zs{t})_\zs{t\ge 0}$ is a standard Brownian motion
and, therefore, $(\dot{w}_\zs{t})_\zs{t\ge 0}$ is ``white noise''. 
By putting here $t=u/\ve$ and $x^{\ve}(t)=y^{\ve}(t\ve)$ we obtain
on the time interval $[0,1]$ the following stochastic differential equation with small
parameter
$$
 \ddot{x}^{\ve}(t)
=
\ve^2
(
p x^{\ve}(t)
+f(t)
)
+
\ve^2
\delta(t) \, \dot{W}_\zs{t}\,,
$$
where $W_\zs{t}=w_\zs{t\ve}/\sqrt{\ve}$ is a standard Brownian motion as well.
In this paper we consider a more general nonlinear case, i.e. the equation of 
the following form
\begin{equation}\label{sec:2.1}
\d \dot{x}^{\ve}(t)
=
\ve^2
(
A(t,x^{\ve}(t))
+f(t)
)
\,\d t
+
\ve^2
\delta (t)\, \d W_\zs{t}\,,
\end{equation}
where $A$ is some $\bbr_\zs{+}\times\bbr\to\bbr$ nonrandom function and 
 $(W_\zs{t})_\zs{0\le t\le 1}$ is the standard Brownian motion.

The problem is to study the asymptotic (as $\ve\to 0$) behavior of the equation
 \eqref{sec:2.1}.
First of all we have to provide an existence
and unicity theorem for strong solutions for nonlinear stochastic equations
of the second order with the
boundary conditions \eqref{sec:2.2} -- \eqref{sec:2.4}. To this end we make use of a 
some uniform version of 
the implicit function theorem. Moreover, to study the asymptotic behavior we propose 
an averaging method based on the 
Green functions approach introduced in  \cite{KaNiQu} for deterministic
differential equations.

The paper is organized as follows. In Section~\ref{sec:3}  we state existence and unicity
theorem. In Section~\ref{sec:Av} we state averaging theorems. 
In Section~\ref{sec:5} we study the corresponding
Green functions. In section~\ref{sec:6} we give the principal proofs. 
In Appendix we prove some technical results.

\section{Existence and unicity theorem}\label{sec:3}

In this section we consider the stochastic equation \eqref{sec:2.1}
with an arbitrary Brownian motion. The first problem is to find sufficient conditions
for existence of unique strong solution. To this end we assume that the function $A(t,x)$ satisfies the following conditions.

\noindent $\C_\zs{1})$ {\sl There exists a positive constant $p$ such that
\begin{equation}\label{sec:3.1}
A(t,x)=px+B(t,x)\,,
\end{equation}
where the function $B(t,x)$ is bounded, i.e.
$$
\beta^{*}=
\sup_\zs{0\le t\le 1}\sup_\zs{x\in\bbr}
|B(t,x)|\,<\,\infty\,.
$$
}

\noindent $\C_\zs{2})$ {\sl
The partial derivative $B^{\prime}_\zs{x}(t,x)$ is bounded by $p$, i.e.
\begin{equation}\label{sec:3.2}
\beta=
\sup_\zs{0\le t\le 1}\,\sup_\zs{x\in\bbr}
\left|
B^{\prime}_\zs{x}(t,x)
\right|<p\,,
\end{equation}
and, moreover, it is uniformly continuous, i.e.
\begin{equation}\label{sec:3.3}
\lim_\zs{\delta\to 0}\,
\sup_\zs{0\le t\le 1}\,
\sup_\zs{|x-y|\le \delta}
\left|
B^{\prime}_\zs{x}(t,x)
-
B^{\prime}_\zs{x}(t,y)
\right|\,=\,0\,.
\end{equation}
}

\noindent For example, we can take $B(t,x)=\beta^{*} \sin(\omega t x)$ for  
$0\le \beta^{*} |\omega| < p$.

\begin{theorem}\label{Th.sec:3.1}
Assume that the conditions $\C_\zs{1})$--$\C_\zs{2})$ hold.
Then there exists nonrandom parameter $\ve_\zs{0}>0$ such that for all $0<\ve\le \ve_\zs{0}$
the equation \eqref{sec:2.1} 
with one of the boundary conditions \eqref{sec:2.2}, \eqref{sec:2.3} or 
 \eqref{sec:2.4} has unique strong continuously differentiable almost sure solution.
\end{theorem}
\noindent The proof of this theorem is given in Section~\ref{sec:6}.

\begin{remark}\label{Re.sec:3.1}
 It should be noted that the well-known example of a stochastic process satisfying 
a stochastic differential equation of a special type with the boundary conditions
\eqref{sec:2.2}
 is the  Brownian Bridge (see, for example, \cite{KaSh}, p. 360), 
which is defined as
\begin{equation}\label{sec:3.4}
Z_\zs{t}=W_\zs{t}-tW_\zs{1}=\int^{1}_\zs{0}G(t,s)\d W_\zs{s}\,,
\end{equation}
where
$$
G(t,s)=(1-t)\Chi_\zs{\{s\le t\}}-t\Chi_\zs{\{s> t\}}\,.
$$
It easy to check directly that this process satisfies the following stochastic equation.
$$
\d Z_\zs{t}=t^{-1}\left(
Z_\zs{t}-W_\zs{t}
\right)\d t
+
\d W_\zs{t}\,,
\quad Z_\zs{0}=Z_\zs{1}=0\,.
$$

\noindent 
One of the possible extensions of the Brownian bridge 
is the backward stochastic equations. (see, for example, in in \cite{ElMa}).
In this paper we 
extent the usual and backward stochastic equations
by introducing 
the
stochastic differential equations of the second order of the form \eqref{sec:2.1}
 with the boundary conditions.

\end{remark}

\begin{remark}\label{Re.sec:3.2}
Note also that we can use the equation \eqref{sec:2.1} - \eqref{sec:1.0}
in the bond markets (see, for example, Chapter 6 in \cite{LamLap}) for the modeling 
of the risk asset as
\begin{equation}\label{sec:3.5}
\B^{\ve}_\zs{t}=e^{(1-t)\ln \B_\zs{0}+x^{\ve}_\zs{t}}\,,
\end{equation}
where $\B_\zs{0}>0$ is any fixed initial price.
\end{remark}

\bigskip

\section{Averaging theorems}\label{sec:Av}

In this section we study the asymptotic (as $\ve\to 0$) properties of the solutions
 of the equation \eqref{sec:2.1} with boundary conditions \eqref{sec:2.2}--\eqref{sec:2.4}.
To state the first theorem we set the following process
\begin{equation}\label{sec:Av.1}
\varkappa(t)=
\int^{1}_\zs{0}\,\Upsilon(t,s)\left(B(s,0)+f(s)\right)\d s
+
\int^1_\zs{0}\,
\Upsilon(t,s)\delta(s)\,\d W_\zs{s}\,,
\end{equation}
where $\Upsilon(t,s)=-
\min(t,s)(1-\max(t,s))$.

 In the sequel we denote by $|\cdot|_\zs{*}$
the uniform norm in $\C([0,1])$ and in  
$\C\left([0,1]\times [0,1]\right)$, i.e. for any $x\in\C([0,1])$
and $G\in \C\left([0,1]\times [0,1]\right)$
\begin{equation}\label{sec:Av.1-0}
|x|_\zs{*}=\sup_\zs{0\le t\le 1}\,|x(t)|\quad\mbox{and}\quad
|G|_\zs{*}=\sup_\zs{0\le t,s\le 1}\,|G(t,s)|\,.
\end{equation}

\noindent 
Moreover, for any function $G$ from $\C([0,1]\times [0,1])$ we set
\begin{equation}\label{sec:Av.1-1}
|G|_\zs{1,*}=\sup_\zs{s\ne t, s,t\in [0,1]}\,|G(t,s)|\,.
\end{equation}

\begin{theorem}\label{Th.sec:Av.1}
Assume that the conditions $\C_\zs{1})$--$\C_\zs{2})$ hold.
Then the solution of the problem \eqref{sec:2.1}--\eqref{sec:2.2}
possess the following limiting form
$$
\P-
\lim_\zs{\ve\to 0}
\,
\left|
\ve^{-2}\,x^{\ve}
-
\varkappa
\right|_\zs{*}\,=\,0\,.
$$
\end{theorem}

\medskip

\noindent Now we introduce the function

\begin{equation}\label{sec:Av.2}
B_\zs{0}(x)=x+\frac{1}{p}\,\int^1_\zs{0}\,B(s,x)\,\d s
\end{equation}
and we define the random variable (if it exists) 
\begin{equation}\label{sec:Av.3}
\zeta=B^{-1}_\zs{0}(\eta)\,,
\end{equation}
where
$$
\eta=-\frac{1}{p}\int^1_\zs{0}\,f(s)\,\d s\,-\,
\frac{1}{p}\,\int^{1}_\zs{0}\,\delta(s)\,\d W_\zs{s}\,.
$$

\begin{theorem}\label{Th.sec:Av.2}
Assume that the conditions $\C_\zs{1})$--$\C_\zs{2})$ hold.
Then the the function \eqref{sec:Av.2} is invertible on $\bbr$. Moreover,
the solutions of the problems \eqref{sec:2.1}--\eqref{sec:2.3}
and \eqref{sec:2.1}--\eqref{sec:2.4} satisfy the following property
$$
\P-
\lim_\zs{\ve\to 0}
\,
\left|
x^{\ve}
-
\zeta
\right|_\zs{*}\,=\,0\,.
$$
\end{theorem}
\begin{remark}\label{Re.sec:Av.0}
It is easy to see, that for $\delta\equiv 0$ and 
$B(t,x) \equiv 0$ we obtain the convergence
\eqref{sec:2.4-1}.
\end{remark}
\begin{remark}\label{Re.sec:Av.1}
It should be noted that the theorems~\ref{Th.sec:3.1}, 
~\ref{Th.sec:Av.1}-- \ref{Th.sec:Av.2} are true for any finite dimension state space also.
\end{remark}

\medskip

\medskip

\section{Properties of the Green functions}\label{sec:5}

Let $p>0$ be arbitrary fixed constant. We make use of the following
differential equation
\begin{equation}\label{sec:5.1}
\ddot{u}-\ve^2 pu=0\,.
\end{equation}

\noindent One can check directly, that for the problem \eqref{sec:5.1} -- \eqref{sec:2.2} 
the Green function is defined as

\begin{equation}\label{sec:5.2}
G_\zs{1,\ve}(t,s)=-
\frac{e^{-\ve\sqrt{p}|t-s|}
g_\zs{1,\ve}(t,s)\,g_\zs{2,\ve}(t,s)}{2\ve\sqrt{p}(1-e^{-2\ve\sqrt{p}})}
\end{equation}
where
$$
g_\zs{1,\ve}(t,s)=1-e^{-2\ve\sqrt{p}\min(t,s)}
\quad\mbox{and}\quad
g_\zs{2,\ve}(t,s)=1-e^{-2\ve\sqrt{p}(1-\max(t,s))}\,.
$$

\begin{proposition}\label{Pr.sec:Av.1}
The function \eqref{sec:5.2} satisfies the following limit properties
\begin{equation}\label{sec:5.3}
\lim_\zs{\ve\to 0}
\,\left|G_\zs{1,\ve}-\Upsilon \right|_\zs{*}\,=\,0
\end{equation}
and
\begin{equation}\label{sec:5.3-0}
\lim_\zs{\ve\to 0}
\,\left|
\frac{\partial }{\partial t}\,
\left(
G_\zs{1,\ve} 
-
\,\Upsilon 
\right)
\right|_\zs{1,*}\,=\,0\,,
\end{equation}
where the norms 
$|\cdot|_\zs{*}$ and $|\cdot|_\zs{1,*}$ are defined in \eqref{sec:Av.1-0} and
\eqref{sec:Av.1-1}.
\end{proposition}

\noindent For the problem \eqref{sec:5.1} -- \eqref{sec:2.3} the Green function is

\begin{equation}\label{sec:5.4}
G_\zs{2,\ve}(t,s)=-
\frac{e^{-\ve\sqrt{p}|t-s|}
g_\zs{3,\ve}(t,s)\,g_\zs{4,\ve}(t,s)}{2\ve\sqrt{p}(1-e^{-2\ve\sqrt{p}})}\,,
\end{equation}
where
$$
g_\zs{3,\ve}(t,s)=1+e^{-2\ve\sqrt{p}\min(t,s)}
\quad\mbox{and}\quad
g_\zs{4,\ve}(t,s)=1+e^{-2\ve\sqrt{p}(\max(t,s))}\,.
$$

\noindent Moreover, for the problem \eqref{sec:5.1} -- \eqref{sec:2.4}
the Green function has the form

\begin{equation}\label{sec:5.5}
G_\zs{3,\ve}(t,s)=-
\frac{e^{-\ve\sqrt{p}|t-s|}
+
e^{-\ve\sqrt{p}(1-|t-s|)}
}{2\ve\sqrt{p}(1-e^{-\ve\sqrt{p}})}\,.
\end{equation}

\begin{proposition}\label{Pr.sec:Av.2}
The functions \eqref{sec:5.4} and  \eqref{sec:5.5} 
satisfy the following limit properties
\begin{equation}\label{sec:5.6}
\lim_\zs{\ve\to 0}
\max_\zs{2\le i\le 3}
\,\left|\ve^2\,G_\zs{i}+
\frac{1}{p}
\right|_\zs{*}\,=\,0
\end{equation}
and
\begin{equation}\label{sec:5.7}
\lim_\zs{\ve\to 0}\,
\ve^2\,
\max_\zs{2\le i\le 3}
\,\left|
\frac{\partial }{\partial t}\,G_\zs{i,\ve}
\right|_\zs{1,*}\,=\,0\,.
\end{equation}
\end{proposition}

\section{Proofs}\label{sec:6}

\subsection{Proof of Theorem~\ref{Th.sec:3.1}}
We define for $1\le i\le 3$ the following $\cC[0,1]\to \cC[0,1]$ operators
\begin{equation}\label{sec:6.1}
\Psi_\zs{i}(\ve,x)(t)=
x(t)
-
\int^1_\zs{0}\,
\wt{G}_\zs{i,\ve}(t,s)\,
B(s,x(s))\,\d s\,,
\end{equation}
where $\wt{G}_\zs{i,\ve}(t,s)=\ve^2 G_\zs{i,\ve}(t,s)$ and
 the corresponding Green functions $G_\zs{i,\ve}$ are defined in 
\eqref{sec:5.2}, \eqref{sec:5.3} and \eqref{sec:5.4} respectively and the function 
$B(t,x)$ is given in the condition $\C_\zs{1})$.

To show this theorem we have to check the conditions of Theorem~\ref{Th.sec:A.1}.
Note that 
$$
\Psi_\zs{1}(0,x)(t)=x(t)
$$
and
\begin{equation}\label{sec:6.2}
\Psi_\zs{2}(0,x)(t)=\Psi_\zs{3}(0,x)(t)=
x(t)+\frac{1}{p}\,\int^1_\zs{0}\,B(s,x(s))\,\d s\,.
\end{equation}
Therefore, the conditions $\AP_\zs{1})$--$\AP_\zs{2})$ are obvious for $i=1$.
Let us check these conditions for $i=2$ and $i=3$.
To this end we introduce for any function $\phi\in\cC[0,1]$ the $\bbr\to\bbr$
function
$$
T_\zs{\phi}(v)=-\frac{1}{p}\,\int^1_\zs{0}\,B(s,\phi(s)+v)\,\d s\,.
$$
It is clear that
$$
T^{\prime}_\zs{\phi}(v)=
-
\frac{1}{p}\,\int^1_\zs{0}\,B^{\prime}_\zs{x}(s,\phi(s)+v)\,\d s
$$
and by the inequality \eqref{sec:3.2}
$$
\sup_\zs{\phi\in \cC[0,1]}\,
\sup_\zs{v\in\bbr}
|T^{\prime}_\zs{\phi}(v)|\,<\,1\,.
$$
Therefore, for any $\phi\in \cC[0,1]$ the equation
\begin{equation}\label{sec:6.3}
v=T_\zs{\phi}(v)
\end{equation}
has an unique solution. This implies directly that in tis case for any $\phi\in\cC[0,1]$
the inverse function in the condition $\AP_\zs{1})$ is given as
$$
x^{\phi}_\zs{0}(t)=\phi(t)+v_\zs{\phi}\,,
$$
where $v_\zs{\phi}$ is the solution of the equation \eqref{sec:6.3}.

As to the condition $\AP_\zs{2})$, note that, the Fr\'echet derivative of
the function \eqref{sec:6.2} is
given by the following linear $\cC[0,1]\to\cC[0,1]$ operator
$$
D_\zs{x}(h)(t)=h(t)+\frac{1}{p}\int^1_\zs{0}\,
B^{\prime}_\zs{x}(s,x(s))h(s)\,\d s\,.
$$

Taking into account the inequality \eqref{sec:3.1} we can directly check that this
operator is isomorphism and
$$
D^{-1}_\zs{x}(h)(t)=h(t)
-
\frac{\int^1_\zs{0}B^{\prime}_\zs{x}(s,x(s))h(s)\,\d s }{p+\int^1_\zs{0}B^{\prime}_\zs{x}(s,x(s))\,\d s}\,.
$$
Now we calculate directly, that for any $x\in\C[0,1]$
$$
\left\|
D^{-1}_\zs{x}
\right\|
=
\sup_\zs{h\in\C[0,1]\,,\,|h|_\zs{*}=1}
\left|
D^{-1}_\zs{x}(h)
\right|
\,
\le\,
1+
\frac{\beta}{p-\beta}
\,.
$$
This implies condition $\AP_\zs{2})$.

Moreover, the boundedness of the function $B(t,x)$ given in the condition $\C_\zs{1})$
implies the condition $\AP_\zs{3})$ for all $1\le i\le 3$. Condition $\C_\zs{2})$ enables
the conditions $\AP_\zs{4})$ and $\AP_\zs{5})$. Therefore, by making use of
Theorem~\ref{Th.sec:A.1} for the random functions
$$
\phi_\zs{i,\ve}(t)=\int^1_\zs{0}\wt{G}_\zs{i,\ve}(t,s)f(s)\d s+
\int^1_\zs{0}\wt{G}_\zs{i,\ve}(t,s)\,\delta(s)\, \d W_\zs{s}
$$
we obtain that there exists some nonrandom parameter 
$\ve_\zs{*}>0$ such that for any $0\le \ve\le \ve_\zs{*}$ and
 for any $1\le i\le 3$ there exist random functions
$x^{\ve}_\zs{i}\in\cC[0,1]$ for which
$$
\Psi_\zs{i,\ve}(\ve,x^{\ve}_\zs{i})=\phi_\zs{i,\ve}\,,
$$
i.e.
\begin{align}\nonumber
x^{\ve}_\zs{i}(t)&=\int^1_\zs{0}\wt{G}_\zs{i,\ve}(t,s)B(s,x^{\ve}_\zs{i}(s))\,\d s\,+\,
\int^1_\zs{0}\wt{G}_\zs{i,\ve}(t,s)f(s)\d s\\[2mm] \label{sec:6.4}
&+
 \int^1_\zs{0}\wt{G}_\zs{i,\ve}(t,s)\,\delta(s)\,\d W_\zs{s}\,.
\end{align}

\noindent 
By Lemma~\ref{Le.sec:A.2} this function satisfies the stochastic differential equation
\eqref{sec:2.1} with $i$th boundary condition in \eqref{sec:2.2}--\eqref{sec:2.4}.

Now we show that this solution is unique. To this end note that through
Lemma~\ref{Le.sec:A.2} we obtain that the uniqueness of the equation \eqref{sec:2.1} is equivalent
 to the uniqueness of the equation \eqref{sec:6.4}. Moreover, taking into account
the asymptotic properties \eqref{sec:5.3}, \eqref{sec:5.6} and
the inequality \eqref{sec:3.2} we can find some parameters $0<\ve_\zs{0}\le \ve_\zs{*}$
and $0<\theta<1$ such that 
$$
\max_\zs{1\le i\le 3}\,
\sup_\zs{0<\ve\le \ve_\zs{0}}\,
\sup_\zs{0\le t,s\le 1}\,
\sup_\zs{x,y\in\bbr}\,
|\wt{G}_\zs{i,\ve}(t,s)|\,
\frac{|B(s,x)-B(s,y)|}{|x-y|}
\,\le \,\theta\,.
$$
This implies immediately that the equation \eqref{sec:6.4} has unique solution.
Hence Theorem~\ref{Th.sec:3.1}.
\endproof
\medskip
\subsection{Proof of Theorem~\ref{Th.sec:Av.1}}

First, we set
$$
\overline{x}^{\ve}(t)=\frac{1}{\ve^{2}}
x^{\ve}_\zs{1}(t)
-
\varkappa(t)
\quad\mbox{and}\quad
\overline{G}_\zs{\ve}(t,s)=
G_\zs{1,\ve}(t,s)
-
\Upsilon(t,s)\,.
$$
By making use of the representation \eqref{sec:6.4} for $i=1$, we obtain 
\begin{equation}\label{sec:6.4-0}
\overline{x}^{\ve}(t)
=\int^1_\zs{0}\left(G_\zs{1,\ve}(t,s)B(s,x^{\ve}_\zs{1}(s))
-
\Upsilon(t,s)B(s,0)
\right)\,\d s\,
+\varsigma^{\ve}_\zs{t}\,,
\end{equation}
where $\varsigma^{\ve}_\zs{t}=\int^1_\zs{0}\,\overline{G}_\zs{1,\ve}(t,s)f(s)\d s+
 \eta^{\ve}_\zs{t}$ and
$$
\eta^{\ve}_\zs{t}=\int^{1}_\zs{0}\,\overline{G}_\zs{\ve}(t,s)\,\delta(s)\,\d W_\zs{s}\,.
$$
It is clear that the property \eqref{sec:5.3}
yields 
$$
\lim_\zs{\ve\to 0}\,
\sup_\zs{0\le t\le 1}\,
\int^1_\zs{0}\,\overline{G}_\zs{1,\ve}(t,s)f(s)\d s
=0\,.
$$
Furthermore, taking into account that 
$$
\Delta \overline{G}_\zs{\ve}(t,t)
=\Delta \overline{G}_\zs{\ve}(t,t)
-
\Delta \overline{G}_\zs{\ve}(t-0,t)
=0\,,
$$
\noindent Lemma~\ref{Le.sec:A.1} yields
$$
\eta^{\ve}_\zs{t}=\eta^{\ve}_\zs{0}
+
\int^{t}_\zs{0}\,\overline{D}_\zs{\ve}(s)\d s\,,
$$
where
$$
\overline{D}_\zs{\ve}(s)
=
\int^1_\zs{0}\, \overline{G}_\zs{1,\ve}(t,s)\,\delta(s)\,\d W_\zs{s}
\quad\mbox{and}\quad
\overline{G}_\zs{1,\ve}(t,s)=
\frac{\partial }{\partial t}\,\overline{G}_\zs{\ve}(t,s)\,.
$$
First of all, note that in view of the property
\eqref{sec:5.3}
$$
\E\,\left(\eta^{\ve}_\zs{0}\right)^{2}
=
\int^{1}_\zs{0}\,\overline{G}^{2}_\zs{\ve}(0,s)\,\delta^{2}(s)\,\d s
\le\,
|\overline{G}_\zs{\ve}|^{2}_\zs{*}\,
\|\delta\|^{2}
\,\to\,0
$$
as $\ve\to 0$. Moreover, 
$$
\sup_\zs{0\le s\le 1}
\E\left(\overline{D}_\zs{\ve}(s)\right)^{2}
\le
\, 
\left|
\overline{G}_\zs{1,\ve}\right|^{2}_\zs{1,*}
$$
and, therefore, in view of the property \eqref{sec:5.3-0} 
$$
\lim_\zs{\ve\to 0}\,\sup_\zs{0\le s\le 1}
\E\left(\overline{D}_\zs{\ve}(s)\right)^{2}=0\,.
$$
Thus,
$$
\P-\lim_\zs{\ve\to 0}\,
|\eta^{\ve}|_\zs{*}
=0
\quad\mbox{and}\quad
\P-\lim_\zs{\ve\to 0}\,
|\varsigma^{\ve}|_\zs{*}
=0\,.
$$
Now we rewrite the equality \eqref{sec:6.4-0} as
\begin{align}\nonumber
\overline{x}^{\ve}(t)
&=\int^1_\zs{0}
\,
\Upsilon(t,s)\,
\left(
B(s,x^{\ve}_\zs{1}(s))
-
B(s,0)
\right)\,\d s\,\\[2mm]\label{sec:6.4-1}
&+
\int^1_\zs{0}\,
\overline{G}_\zs{\ve}(t,s)B(s,x^{\ve}_\zs{1}(s))
\,\d s
+\varsigma^{\ve}_\zs{t}\,.
\end{align}
Taking into account here  the conditions $\C_\zs{1}$) and $\C_\zs{2}$)
we obtain that
\begin{align*}
|\overline{x}^{\ve}|_\zs{*}&\le 
\,\beta\,|\Upsilon|_\zs{*}\, |x^{\ve}_\zs{1}|_\zs{*}
+\beta^{*}\,|\overline{G}_\zs{\ve}|_\zs{*}
+|\varsigma^{\ve}|_\zs{*}\\[2mm]
&
\le 
\ve^{2}
\beta\,|\Upsilon|_\zs{*}\, |\overline{x}^{\ve}|_\zs{*}
+
\ve^{2}
\beta\,|\Upsilon|_\zs{*}\, |\varkappa|_\zs{*}
+\beta^{*}\,|\overline{G}_\zs{\ve}|_\zs{*}
+|\varsigma^{\ve}|_\zs{*}\,.
\end{align*}

\noindent
Therefore, for sufficiently small $\ve$ we get
$$
|\overline{x}^{\ve}|_\zs{*}\le
\frac{1}{1-\ve^{2}
\beta\,|\Upsilon|_\zs{*}} 
\left(
\ve^{2}
\beta\,|\Upsilon|_\zs{*}\, |\varkappa|_\zs{*}
+\beta^{*}\,|\overline{G}_\zs{\ve}|_\zs{*}
+|\varsigma^{\ve}|_\zs{*}
\right)
\,\to 0
$$
as $\ve\to 0$. Hence Theorem~\ref{Th.sec:Av.1}.
\endproof

\medskip

\subsection{Proof of Theorem~\ref{Th.sec:Av.2}}
First of all note that the boundedness of the function $B(s,x)$ in the condition $\C_\zs{1})$ 
implies 
$$
\lim_\zs{x\to+\infty}\,B_\zs{0}(x)=+\infty
\quad\mbox{and}\quad
\lim_\zs{x\to-\infty}\,B_\zs{0}(x)=-\infty\,.
$$
Moreover, by the inequality \eqref{sec:3.2} we get
\begin{equation}\label{sec:6.5}
\b_\zs{*}=\inf_\zs{x\in\bbr}\,|B^{\prime}_\zs{0}(x)|\,>\,0\,.
\end{equation}
These properties imply that for any $a\in \bbr$ the equation
$$
B_\zs{0}(x)=a
$$
has an unique solution, i.e. the function $B_\zs{0}$ defined in \eqref{sec:Av.2}
is invertible on $\bbr$.
Moreover, the representation \eqref{sec:6.4} implies
\begin{equation}\label{sec:6.6}
B_\zs{0}(x^{\ve}(t))=\eta+
u^{\ve}(t)
+
v^{\ve}(t)\,,
\end{equation}
where
$$
u^{\ve}(t)=\int^1_\zs{0}\,V^{\ve}(t,s)\,
\left(
B(s,x^{\ve}(s))
+
f(s)
\right)\,\d s
$$
and
$$
v^{\ve}(t)=\int^1_\zs{0}\,V^{\ve}(t,s)\,\delta(s)\,
\d W_\zs{s}
$$
with
$$
V^{\ve}(t,s)=\wt{G}(t,s)+\frac{1}{p}\,.
$$

\noindent The limit \eqref{sec:5.6} implies directly
\begin{equation}\label{sec:6.7}
\P-\lim_\zs{\ve\to 0}\,
\sup_\zs{0\le t\le 1}\,
|u^{\ve}(t)|\,=\,0\,.
\end{equation}
Now we show that

\begin{equation}\label{sec:6.8}
\P-\lim_\zs{\ve\to 0}\,
\sup_\zs{0\le t\le 1}\,
|v^{\ve}(t)|\,=\,0\,.
\end{equation}

Indeed, by Lemma~\ref{Le.sec:A.1} process $v^{\ve}(t)$ is almost sure
continuously differentiable and taking into account that 
$\Delta V^{\ve}(t,t)=0$ we find that
$$
v^{\ve}(t)
=
v^{\ve}(0)
+
\int^{t}_\zs{0}\,
\check{D}_\zs{\ve}(s)\,\d s\,,
$$
where
$$
\check{D}_\zs{\ve}(s)=
\int^1_\zs{0}\,
\frac{\partial }{\partial s}\,
V^{\ve}(s,\vartheta)\,\delta(\vartheta)\,\d W_\zs{\vartheta}
=
\int^1_\zs{0}\,
\frac{\partial }{\partial s}\,
\wt{G}^{\ve}(s,\vartheta)\,\delta(\vartheta)\,\d W_\zs{\vartheta}\,.
$$

\noindent Moreover, by \eqref{sec:5.6}
$$
\E \left(
v^{\ve}(0)
\right)^2=
\int^1_\zs{0}\,\left(
V^{\ve}(0,s)
\right)^2\,\d s
\to 0
\quad\mbox{as}\quad\ve\to 0
$$
and by \eqref{sec:5.7}
$$
\sup_\zs{0\le s\le 1}\,
\E\left(
\check{D}_\zs{\ve}(s)
\right)^2
=
\sup_\zs{0\le s\le 1}\,
\int^1_\zs{0}\,
\left(
\frac{\partial }{\partial s}\,
\wt{G}^{\ve}(s,\vartheta)
\right)^2\,\d \vartheta
\to 0
\quad\mbox{as}\quad
\ve\to 0\,.
$$
This implies immediately \eqref{sec:6.8}.
Therefore, taking into account that 
$$
B_\zs{0}(\zeta)=\eta
$$ 
and that
$$
|x-y|\le \frac{1}{\b_\zs{*}}\,|B_\zs{0}(x)-B_\zs{0}(y)|
$$
we obtain Theorem~\ref{Th.sec:Av.2}.
\endproof

\medskip

\section{Conclusion}\label{sec:Co}

In this paper we introduced and studied the stochastic differential equations 
of the second order containing a small parameter with boundary conditions. 
We proved an existence and unicity theorem for strong solutions and we shown 
an average principle as the small parameter goes to zero. Through the Green functions 
method we obtained the explicit limit forms for solutions.

\medskip

\renewcommand{\theequation}{A.\arabic{equation}}
\renewcommand{\thetheorem}{A.\arabic{theorem}}
\renewcommand{\thesubsection}{A.\arabic{subsection}}
\section{Appendix}\label{sec:A}
\setcounter{equation}{0}
\setcounter{theorem}{0}

\subsection{Technical Lemmas}\label{subsec:A.1}
\begin{lemma}\label{Le.sec:A.1}
Let $G(t,s)$ be a function having 
the following form
$$
G(t,s)=
\left\{
\begin{array}{cc}
\sum^m_\zs{j=1}g_\zs{1j}(t)h_\zs{1j}(s)
&\quad\mbox{if}\quad t\ge s\,;\\[4mm]
\sum^m_\zs{j=1}g_\zs{2j}(t)h_\zs{2j}(s)
&\quad\mbox{if}\quad t<s\,.
\end{array}
\right.
$$
Assume that the functions $g_\zs{1j}$ and $g_\zs{2j}$
belong to $\C^2[0,1]$. Then the stochastic process
$$
\xi_\zs{t}=\int^1_\zs{0}G(t,s)\,\delta(s)\,\d W_\zs{s}
$$
admits the following Ito differential 
\begin{equation}\label{sec:A.1}
\d \xi_\zs{t}=
D(t)
\d t
+
\Delta G(t,t)\,\delta(t)\,\d W_\zs{t}\,,  
\end{equation}
where
$$
D(t)
=
\int^1_\zs{0}\, G_\zs{1}(t,s)\,\delta(s)\,\d W_\zs{s}
\quad\mbox{with}\quad
G_\zs{1}(t,s)=
\frac{\partial }{\partial t}\,G(t,s)
$$
and
$$
\Delta G(t,t)=G(t,t)-G(t-0,t)=
\sum^m_\zs{j=1}g_\zs{1j}(t)h_\zs{1j}(t)
-
\sum^m_\zs{j=1}g_\zs{2j}(t)h_\zs{2j}(t)\,.
$$
\end{lemma}
\proof
First, note that
$$
\xi_\zs{t}=\int^{t}_\zs{0}G(t,s)\,\delta(s)\,\d W_\zs{s}
+\int^1_\zs{t}G(t,s)\,\delta(s)\,\d W_\zs{s}\,.
$$
Taking into account here the definition of the function 
$G$ we get
\begin{align*}
 \xi_\zs{t}
&=\sum^m_\zs{j=1}g_\zs{1j}(t)\,\int^{t}_\zs{0}\,h_\zs{1j}(s)\,\d W_\zs{s}
+\sum^m_\zs{j=1}g_\zs{2j}(t)\,\int^{1}_\zs{t}\,h_\zs{2j}(s)\,\d W_\zs{s}\\[2mm]
&=\sum^m_\zs{j=1}g_\zs{1j}(t)\,\int^{t}_\zs{0}\,h_\zs{1j}(s)\,\d W_\zs{s}
+\sum^m_\zs{j=1}g_\zs{2j}(t)\,\int^{1}_\zs{0}\,h_\zs{2j}(s)\,\d W_\zs{s}\\[2mm]
&-\sum^m_\zs{j=1}g_\zs{2j}(t)\,\int^{t}_\zs{0}\,h_\zs{2j}(s)\,\d W_\zs{s}\,.
\end{align*}
Now, this Lemma follows directly from the Ito formula.
\endproof

Now we make use of the following conditions for the Green functions

\noindent $\H_\zs{1})$ {\sl The 
$[0,1]\times[0,1]\to \bbr$ function $G(\cdot,\cdot)$ is continuous and has the following form
\begin{equation}\label{sec:A.2}
G(t,s)=
\left\{
\begin{array}{cc}
g_\zs{11}(t)h_\zs{11}(s)+
g_\zs{12}(t)h_\zs{12}(s)
&\quad\mbox{if}\quad t\ge s\,;\\[2mm]
g_\zs{21}(t)h_\zs{21}(s)+
g_\zs{22}(t)h_\zs{22}(s)
&\quad\mbox{if}\quad t< s\,,
\end{array}
\right.
\end{equation}
where the functions $g_\zs{11}(\cdot)$,  $g_\zs{12}(\cdot)$, $g_\zs{21}(\cdot)$
$g_\zs{22}(\cdot)$ are two times continuously differentiable such that
\begin{equation}\label{sec:A.3}
\frac{\partial^2 }{\partial t^2} G(t,s)-
p(t)\,G(t,s)=0
\end{equation}
for any $t\ne s$ from $[0,1]$.
}

\noindent $\H_\zs{2})$ {\sl The first partial derivative of the function $G$
has the jump
$$
\frac{\partial }{\partial t}\,G(t,t)
-
\frac{\partial }{\partial t}\,G(t-0,t)
=1\,.
$$
}

\begin{lemma}\label{Le.sec:A.2}
Assume that the function $G$ satisfies the conditions 
$\H_\zs{1}$)--$\H_\zs{2}$) and $\varphi$ is arbitrary  function from
$\C[0,1]$. Then the stochastic process
$$
\xi_\zs{t}=
\int^1_\zs{0}\,G(t,s)\,\varphi(s)\,\d s+
\int^1_\zs{0}\,G(t,s)\,\delta(s)\,\d W_\zs{s}
$$
is continuously almost sure differentiable
and 
satisfies the following second order stochastic
equation
\begin{equation}\label{sec:A.4}
\ddot{\xi}_\zs{t}=p(t)\xi_\zs{t}+
\varphi(t)+\delta(t)\,
\dot{W}_\zs{t}\,,
\end{equation}
i.e. the derivative $\dot{\xi}_\zs{t}$ satisfies the following stochastic
equation
$$
\d \dot{\xi}_\zs{t}=
p(t) \xi_\zs{t}\d t\,
+\,
\varphi(t)\,
\d t\,
+\,\delta(t)\d W_\zs{t}
\,.
$$
\end{lemma}
\proof
First of all we set
\begin{equation}\label{sec:A.4-0}
\a_\zs{t}=\int^1_\zs{0}\,G(t,s)\,\varphi(s)\,\d s
\quad\mbox{and}\quad
\b_\zs{t}=
\int^1_\zs{0}\,G(t,s)\,\delta(s)\,\d W_\zs{s}\,.
\end{equation}
It is easy to check directly that
$$
\ddot{\a}_\zs{t}=
p(t)\a_\zs{t}+\varphi(t)\,.
$$
Moreover, note that Condition $\H_\zs{1})$ implies $\Delta G(t,t)=0$. 
Therefore, by Lemma~\ref{Le.sec:A.1} we obtain 
$$
\b_\zs{t}=\int^t_\zs{0}\,
D(s)\,\d s\,,
$$
where the process $D(t)$ is defined in \eqref{sec:A.1}. From here we obtain that
$$
\dot{\b}_\zs{t}\,=\,D(t)=\int^1_\zs{0}\,
G_\zs{1}(t,s)\,\delta(s)\,\d W_\zs{s}
\,.
$$
Therefore, by Lemma~\ref{Le.sec:A.1} we find
$$
\d \dot{\b}_\zs{t}=
D_\zs{1}(t)\,\d t
+
\Delta G_\zs{1}(t,t)\,\delta(t)\,\d W_\zs{t}\,,
$$
where
$$
D_\zs{1}(t)=
\int^1_\zs{0}\,\frac{\partial }{\partial t} G_\zs{1}(t,s)\,\delta(s)\,\d W_\zs{s}
=
\int^1_\zs{0}\,\frac{\partial^2 }{\partial t^2} G(t,s)\,\delta(s)\,\d W_\zs{s}
\,.
$$
Now in view of  the equations \eqref{sec:A.3}
and the definition \eqref{sec:A.4-0} we obtain that 
$$
D_\zs{1}(t)=p(t)\b(t)\,.
$$
Therefore, the  condition $\H_\zs{2})$
implies this lemma. \endproof

\subsection{Green functions}\label{subsec:A.2}

In this section we consider the second order 
$\C^{2}[0,1]\to \C^{2}[0,1]$ linear differential operator
$L$ defines as
$$
L(x)(t)=\ddot{x}(t)+p(t)\dot{x}(t)+q(t)x(t)
$$
with the boundary conditions
\begin{equation}\label{sec:A.5}
\alpha_\zs{1}
x(0)
+
\beta_\zs{1}
\dot{x}(0)
=0
\quad\mbox{and}\quad
\alpha_\zs{2}
x(1)
+
\beta_\zs{2}
\dot{x}(1)
=0\,,
\end{equation}
where $\alpha_\zs{1}$, $\beta_\zs{1}$, $\alpha_\zs{2}$ 
and $\beta_\zs{2}$
are some fixed constants. We will consider also the following mixed conditions
\begin{equation}\label{sec:A.6}
x(0)
=
x(1)
\quad\mbox{and}\quad
\dot{x}(0)
=
\dot{x}(1)\,.
\end{equation}

\noindent The operator $L$ is called {\sl regular} if the equation
\begin{equation}\label{sec:A.7}
L(x)=0
\end{equation}
with the boundary conditions \eqref{sec:A.5} has only trivial solution $x\equiv 0$.
If the operator $L$ is regular, then the problem
$$
L(x)=f
$$
with the boundary conditions \eqref{sec:A.5} may be written as
$$
x(t)=\int^1_\zs{0}\,G(t,s)\,f(s)\,\d s\,.
$$
Here the $[0,1]^2\to\bbr$ function $G$ is called {\sl Green function} generated by the operator
$L$ and the boundary conditions \eqref{sec:A.5}. The function $G$ has the following
properties

\begin{enumerate}
 \item {\sl $G(t,s)$ is continuous with respect to $t$ and $s$.}
\item {\sl For $t\ne s$ the function $G(t,s)$ is two times continuously differentiable
with respect to $t$ and $s$. Moreover, for $t\ne s$
$$
\frac{\partial^2}{\partial t^2}\,G(t,s)
+
p(t)
\frac{\partial }{\partial t }\,G(t,s)
+
q(t) G(t,s)=0
$$
and 
$$
\alpha_\zs{1}
G(0,s)
+
\beta_\zs{1}
\frac{\partial }{\partial t}
G(0,s)
=0
\quad\mbox{and}\quad
\alpha_\zs{1}
G(1,s)
+
\beta_\zs{1}
\frac{\partial }{\partial t}
G(1,s)
=0
\,,
$$
or
$$
G(0,s)
=
G(1,s)
\quad\mbox{and}\quad
\frac{\partial }{\partial t}
G(0,s)
=\frac{\partial }{\partial t}
G(1,s)\,.
$$
}

\item
{\sl The function $G$ is symmetric, i.e.}
$$
G(t,s)=G(s,t)
$$
\item {\sl The partial derivative with respect to $t$ has the unitary jump, i.e.
$$
\frac{\partial }{\partial t}
G(t+0,t)
-
\frac{\partial }{\partial t}
G(t-0,t)
=
1\,.
$$
}
\end{enumerate}

\subsection{Uniform implicit function theorem}

Let $F$ be a Banch space and $\Psi\,:\,[0,1]\times F\to F$ be a continuous function.
In this section we study the following equation
\begin{equation}\label{sec:A.8}
\Psi(\ve,x)=\phi
\quad\mbox{for}\quad
\phi\in F\,.
\end{equation}

\noindent $\AP_\zs{1})$
{\sl The equation \eqref{sec:A.8} for $\ve=0$ has unique solution $x^{\phi}_\zs{0}$ for any 
$\phi\in F$.}

\noindent $\AP_\zs{2})$
{\sl For any $\phi\in F$ the partial Fr\'echet derivative
$$
\Upsilon(\phi)=\frac{\partial }{\partial x}\,\Psi(0,x^{\phi}_\zs{0})
$$
is $F\to F$ isomorphism and

\begin{equation}\label{sec:A.9}
\sup_\zs{\phi\in F}\,
\left|
\Upsilon^{-1}(\phi)\,
\right|\,<\,\infty\,.
\end{equation}
}
\noindent $\AP_\zs{3})$
{\sl The function $\Psi$ is uniform continuous in $\ve$, i.e.}
$$
\lim_\zs{\ve\to 0}\,\sup_\zs{x\in F}\,
|\Psi(\ve,x)-\Psi(0,x)|\,=\,0\,.
$$
\noindent $\AP_\zs{4})$
{\sl The function $\Psi$ is uniformly in the Fr\'echet since differentiable with respect to $x$, i.e.
$$
\lim_\zs{\delta\to 0}\,
\sup_\zs{0<\ve\le 1}
\sup_\zs{|x-y|\le \delta}\,
\frac{
\left|
\Psi(\ve,y)
-
\Psi(\ve,x)
-
\dfrac{\partial }{\partial x}\Psi(\ve,x)\,(y-x)
\right|
}
{|x-y|}
=0\,.
$$
}

\noindent $\AP_\zs{5})$
{\sl The partial Fr\'echet derivative of $\partial\Psi/\partial x$ is uniform
continuous, i.e.
$$
\lim_\zs{\ve\to 0}\,
\limsup_\zs{\delta\to 0}\,
\sup_\zs{|x-y|<\delta}
\left|
\frac{\partial }{\partial x}\,\Psi(\ve,x)
-
\frac{\partial }{\partial x}
\,\Psi(0,y)
\right|\,
=\,0\,.
$$
}

\begin{theorem}\label{Th.sec:A.1}
Assume that the conditions $\AP_\zs{1})$--$\AP_\zs{5})$ hold. 
Then  there exist
$\ve_\zs{0}\in (0,1)$
and $\delta>0$ such that for any $0<\ve\le\ve_\zs{0}$
and any $\phi\in F$ there exists unique 
 $x^{\phi}_\zs{\ve}\in \B(x^{\phi}_\zs{0},\delta)$ which satisfies the equation
\eqref{sec:A.8}, i.e.
$$
\Psi(\ve, x^{\phi}_\zs{\ve})=\phi\,,
$$
where $\B(a,\delta)=\{x\in F\,:\,|x-a|\le \delta\}$.
\end{theorem}

\proof
By the same way as in \cite{Ca}(page 57) we introduce the following special 
function
$$
F(\ve,z)=z-\Upsilon^{-1}(\phi)\left(
\Psi(\ve,z)-\phi
\right)\,.
$$
Note that the conditions $\AP_\zs{2})$, $\AP_\zs{4})$ and $\H_\zs{5})$ imply
that there exist $\ve_\zs{*}>0$, $\delta>0$ and $0<\theta<1$ 
such that for any $0\le\ve\le  \ve_\zs{0}$ and $x$, $y$ from $\B(x^{\phi}_\zs{0},\delta)$
$$
|F(\ve,x)-F(\ve,y)|\le \theta |x-y|\,.
$$
Now we define the following approximating sequence
$$
z^{\ve}_\zs{n}=F(\ve,z^{\ve}_\zs{n-1})
\quad\mbox{for}\quad n\ge 1
$$
and $z^{\ve}_\zs{0}=x^{\phi}_\zs{0}$.
First of all note that 
$$
|z_\zs{1}-x^{\phi}_\zs{0}|=|z_\zs{1}-z_\zs{0}|
=
\left|
\Upsilon^{-1}(\phi)
\left(
\Psi(\ve,x^{\phi}_\zs{0})
-
\Psi(0,x^{\phi}_\zs{0})
\right) 
\right|
$$
and by the inequality \eqref{sec:A.9} and the condition $\AP_\zs{3})$ we can choose
$0<\ve_\zs{0}\le \ve_\zs{*}$ such that for all $0<\ve\le\ve_\zs{0}$
$$
|
z^{\ve}_\zs{1}
-
z^{\ve}_\zs{0}
|
\le 
(1-\theta)\delta
<\delta\,,
$$
i.e. $z^{\ve}_\zs{1}\in \B(z_\zs{0},\delta)$. 
Now through the induction method we assume that 
$z^{\ve}_\zs{j}\in \B(z_\zs{0},\delta)$ for all 
$1\le j\le n-1$ with $n\ge 2$. Let us check that in this case 
$z^{\ve}_\zs{n}\in \B(z_\zs{0},\delta)$ as well. Indeed, taking into account
 our condition implies directly the following inequality
$$
|z^{\ve}_\zs{n}-z^{\ve}_\zs{n-1}|
\le \theta 
|z^{\ve}_\zs{n-1}-z^{\ve}_\zs{n-2}|
\le \theta^{n-1}\,|z^{\ve}_\zs{1}-z^{\ve}_\zs{0}|
$$
we get
$$
|z^{\ve}_\zs{n}-x^{\phi}_\zs{0}|
=
|z^{\ve}_\zs{n}-z^{\ve}_\zs{0}|
\le \sum^n_\zs{j=1}\,| z^{\ve}_\zs{j}-z^{\ve}_\zs{j-1}|
\le 
\sum^n_\zs{j=1}\,\theta^{j-1}
|z^{\ve}_\zs{1}-z^{\ve}_\zs{0}|
\le \delta\,.
$$
Thus $z^{\ve}_\zs{n}\in \B(z_\zs{0},\delta)$ for all $n\ge 2$.
Therefore, for any $n\ge 1$
$$
|z^{\ve}_\zs{n}-z^{\ve}_\zs{n-1}|
\le \theta 
|z^{\ve}_\zs{n-1}-z^{\ve}_\zs{n-2}|\,,
$$
i.e. for any $0\le \ve\le \ve_\zs{0}$ there exists the limit
$$
\lim_\zs{n\to\infty}\,z^{\ve}_\zs{n}=z^{\ve}_\zs{*}\,,
$$
which is the unique solution of the equation $z=F(\ve,z)$ in the ball 
$\B(z_\zs{0},\delta)$ for $\ve\le \ve_\zs{0}$. This implies immediately
Theorem~\ref{Th.sec:A.1}. \endproof


\medskip

\newpage
\begin{thebibliography}{100}

\bibitem{AlBa}
Alvarez, O. and Bardi, M. (2001/02)
Viscosity solutions methods for singular perturbations in deterministic 
and stochastic control. {\sl SIAM J. Control Optim.} {\bf 40} (4), 
1159–1188 

\bibitem{ArGa}
Artzein, Z. and Gaitsgory, V. (2000)
The value function of singularly perturbed control systems.
{\sl Applied Mathematics and Optimization}, {\bf 41}, 425 -- 445.

\bibitem{Be}
Bensoussan, A. {\sl Perturbation Methods in Optimal Control.}
J. Wiley/Gauthier Villars, New York, 1988.

\bibitem{AiMyZh}
Ait-Sahalia, Y., Mykland, P.A. and Zhang, L. (2005)
How often to sample a continuous-time process in the presence of market
microstructure noise. {\sl Rev. Financ. Studies}, {\bf 18}, 351--416.

\bibitem{AiMyZh1}
Ait-Sahalia, Y., Mykland, P.A. and Zhang, L. (2005)
A tale of two times scales: Determining integrated volatility with noisy high-frequency
data. {\sl J. Amer. Statist. Assoc.}, {\bf 100}, 1394--1411.

\bibitem{BoMi}
Bogolubov, N.N. and Mitropol'skii, Yu.A.
{\sl Asymptotic methods in the theory of nonlinear oscillation.}
M.: Fizmatgiz, 1963 (in Russian).




\bibitem{Ca}
Cartan, H. {\sl Cours de calcul diff\'erentiel.}
Hermann, \'Editeurs des sciences et des arts, 1990.

\bibitem{ElMa}
El Karoui, N. and Mazliak, L. {\sl Backward Stochastic differential
equations.} Pitman Research Notes in Mathematics Series, 364. Longman, Harlow, 1997.

\bibitem{FoPaSi}
 Fouque, J-P.,  Papanicolaou, G. and Sicar, K.
{\sl Derivatives in financial markets with stochastic volatility.}
Cambridge University Press, 2000.

\bibitem{FrWe}
Freidlin, M.I. and Wentzell, A.D. 
{\sl Random Perturbations of Dynamical Systems.}
Springer-Verlag, New York, 1986.





\bibitem{KaPe}
Kabanov, Yu. and Pergamenshchikov, S. 
{\sl Two Scale Stochastic 
Systems: Asymptotic Analysis and Control.} Applications 
of Mathematics, Stochastic Modeling and Applied Probability, 49, 
Springer-Verlag, Berlin, New York, 2003.

\bibitem{KaNiQu}
Kamenskii, M., Nistri, P. and Quincampoix, M. (2002)
 Sliding mode control of uncertain systems: a singular perturbation approach. 
{\sl  IMA J. Math. Control Inform.},  {\bf 19} (4), 377--398. 

\bibitem{KaSh}
Karatzas, I. and Shreve, S.E. 
{\sl Brownian Motion and Stochastic Calculus}.
Springer-Verlag, 1991.

\bibitem{KhaKr}
Khasminskii, R. and  Krylov, N. (2001)
On averaging principle for diffusion processes with nullrecurrent
fast component. {\sl Stoch Proc Appl.} {\bf 93},  229-240.

\bibitem{Ku}
Kushner, H. 
{\sl Weak Convergence Methods and Singularly Perturbed
Stochastic Control and Filtering.}
IEEE Press, New York, 1986.

\bibitem{LamLap}
Lamberton, D. and Lapeyre, B.
{\sl Introduction to stochastic calculus applied to finance.}
Chapman \& Hall, London, 1996. 

\bibitem{PaPaSt}
Papavasiliou, A., Pavliotis, G.A., Stuart, A. M. (2009) 
Maximum likelihood drift estimation for multi-scale diffusions. 
{\em Stochastic Process. Appl.}, {\bf 119} (10), 3173--3210. 

\bibitem{PaSt}
Pavliotis, G.A., Stuart, A. M. (2007)
Parameter estimation for multiscale diffusions.
{\sl J. Stat. Phys.}, {\bf 127} (4), 741--781.

\bibitem{PaSt1}
Pavliotis, G.A., Stuart, A. M. (2008)
{\sl Multiscale Methods.}
Texts in Applied Mathematics, Averaging and harmonization {\bf 53}, Springer, New York

\bibitem{Ti}
Tikhonov. A.N. (1950) 
On system of differential equations containing parameters.
{\em Matem. Sbornik}, {\bf 27} (69),  147-156. 

\bibitem{Ti1}
Tikhonov. A.N. (1952) 
System of differential equations containing a small parameters at derivatives.
{\em Matem. Sbornik}, {\bf 31} (73),  147-156. 




\end{thebibliography}
\end{document}